\documentclass[12pt,fleqn,twoside]{article}
\usepackage[cp1251]{inputenc}
\usepackage[ukrainian]{babel}
\usepackage{amsmath,amsfonts,amssymb}
\usepackage{latexsym,hhline,graphics,epsfig}

\newcounter{theorem}
\newcommand{\theor}{\par\refstepcounter{theorem}%
{\bf Теорема \thetheorem }.\,\,}

\newcounter{lemma}
\newcommand{\be}{\begin{equation}}
\newcommand{\ee}{\end{equation}}

\begin{document}

\large \pagestyle{empty}

\noindent {\small \textbf{УДК}\ \ 517.54}

\bigskip


\begin{center}{\bf A.К. Бахтин, Г.П. Бахтина, И.В. Денега}\\\end{center}
\vskip 4mm

(Институт математики НАН Украины, Киев)

\vskip 4mm

(Национальный технический университет Украины "Киевский
политехнический институт", Киев)

\vskip 4mm

(Институт математики НАН Украины, Киев)

\begin{center}{\bf Неравенства в задачах о неналегающих областях.}\end{center}

\vskip 4mm

alexander.bahtin@yandex.ru, iradenega@yandex.ru

\vskip 4mm
\parbox{16.5cm}{\footnotesize В данной работе рассмотрена одна достаточно общая задача о
неналегающих областях со свободными полюсами на лучевых системах.
Основная теорема данной статьи значительно обобщает ранее известные
результаты для задач подобного типа.}

\bigskip
\parbox{16.5cm}{\footnotesize In this paper we consider quite general problem on
non-overlapping domains with free poles on radial systems.
The main theorem of this work generalizes the previously known
results for problems of this type.}

\bigskip

\textbf{1. Введение.} В геометрической теории функций комплексного
переменного экстремальные задачи о неналегающих областях являются
хорошо известным классическим направлением [1 -- 43]. Возникновение
данного направления связано с классической работой М.А. Лаврентьева
\cite{1}, в которой, в частности, была впервые поставлена и решена
задача о максимуме произведения конформных радиусов двух
непересекающихся односвязных областей. В дальнейшем этот результат
получил значительное развитие в работах многих авторов [2 -- 43]. В
последнее время особое внимание специалистов привлекли
задачи с так называемыми "свободными" полюсами \cite{bah}--\cite{38BakhG}, \cite{kyzmina}--\cite{70}. При
исследовании экстремальных задач важную роль играет теория
квадратичных дифференциалов (см.[4]).

В работах \cite{3}--\cite{59}, \cite{10} был, в частности,
разработан метод кусочно-разделяющего преобразования, который открыл
новые возможности для решения экстремальных задач о неналегающих
областях. При определенных условиях этот метод позволяет сводить
задачи с большим числом неизвестных параметров к задачам с меншим их
числом. Данная работа посвящена исследованию некоторых задач
подобного типа.

\newpage
\textbf{2. Обозначения и определения.}

Пусть $\mathbb{N}$, $\mathbb{R}$ -- множества натуральных и
вещественных чисел соответственно, $\mathbb{C}$ -- комплексная
плоскость, $\overline{\mathbb{C}}=\mathbb{C}\bigcup\{\infty\}$ -- ее
одноточечная компактификация,  $r(B,a)$ -- внутрений радиус области
$B\subset\overline{C}$, относительно точки $a\in B$ (см. напр.
\cite{heim}, \cite{5}, \cite{10}), $\chi(t)=\frac{1}{2}(t+t^{-1})$ и
$\mathbb{R^{+}}=(0,\infty)$.

Пусть $n\in \mathbb{N}$, $n\geq 2$. Систему точек
$A_{n}:=\left\{a_{k} \in \mathbb{C}:\, k=\overline{1,n}\right\},$
назовем \textbf{\emph{$n$ - лучевой}}, если
$|a_{k}|\in\mathbb{R^{+}}$ при $k=\overline{1,n}$, $0=\arg
a_{1}<\arg a_{2}<\ldots<\arg a_{n}<2\pi$. Обозначим при этом
$P_{k}(A_{n}):=\{w: \arg a_{k}<\arg w<\arg a_{k+1}\}$,
$\theta_{k}:=\arg a_{k},\,a_{n+1}:=a_{1},\, \theta_{n+1}:=2\pi,$
$\alpha_{k}:=\frac{1}{\pi}\arg \frac{a_{k+1}}{a_{k}},\,
\alpha_{n+1}:=\alpha_{1},\,k=\overline{1, n}.$

Данная работа базируется на применении кусочно-разделяющего
преобразования, развитого в \cite{3}--\cite{59}. Для конкретного
использования этого метода рассмотрим специальную систему конформных
отображений. Пусть
$\zeta=\pi_k(w)=-i\left(e^{-i\theta_k}w\right)^\frac{1}{\alpha_k}$,\quad
$k=\overline{1,n}$ обозначает ту однозначную ветвь многозначной
аналитической функции $\pi_k(w)$, которая осуществляет однолистное и
конформное отображение $P_k(A_{n})$ на правую полуплоскость
$\text{Re}\,\zeta>0$.

Для произвольной $n$-лучевой системы точек $A_{n}=\{a_{k}\}$ и
$\gamma\in\mathbb{R^{+}}$ полагаем
$$\mathcal{L}^{(\gamma)}(A_n):=\prod\limits_{k=1}^n\left[
\chi\left(\Bigl|\frac{a_k}{a_{k+1}}\Bigr|^\frac{1}{2\alpha_k}\right)\right]^{1-\frac{1}{2}\gamma\alpha_k^2}\cdot
\prod\limits_{k=1}^n|a_k|^{1+\frac{1}{4}\gamma(\alpha_k+\alpha_{k-1})}.$$
Класс тех $n$ - лучевых систем точек для которых
$\mathcal{L}^{(\gamma)}(A_n)=1$ автоматически содержит все $n$ -
лучевые системы точек, рассположенные на окружности.

Целью данной работы является получение точных оценок сверху для
функционалов следующего вида
\begin{equation}\label{1a}J_n^{(\gamma)}=r^\gamma\left(B_0,0\right)\prod\limits_{k=1}^n r\left(B_k,a_k\right),\end{equation}
где $\gamma\in\mathbb{R^{+}}$, $A_{n}=\{a_{k}\}_{k=1}^{n}$ --
$n$-лучевая система точек, $a_{0}=0,$ $\{B_{k}\}_{k=0}^{n}$ --
система неналегающих областей (то есть $B_{p}\cap B_{j}={\O}$ при
$p\neq j$, $p,j=\overline{0, n}$) таких, что $a_{k}\in B_{k}\subset \overline{\mathbb{C}}$
при $k=\overline{0, n}$.

Введем функцию $F(x)=2^{x^2+6}\cdot x^{x^2+1}\cdot (2-x)^
{-\frac{1}{2}(2-x)^2} (2+x)^{-\frac{1}{2}(2+x)^2}$, $x \in (0,2]$.

\textbf{3. Основные результаты.} \textbf{\theor}{\it Пусть $n\in
\mathbb{N}$, $n\geqslant 2$ и $\gamma\in (0, 1]$. Тогда для любой
$n$-лучевой системы точек $A_n=\{a_k\}_{k=1}^n$,
$\mathcal{L}^{(\gamma)}\left(A_n\right)=1$ и любого набора взаимно
непересекающихся областей $B_k$, $a_k\in
B_k\subset\overline{\mathbb{C}}$, $k=\overline{0,n}$, справедливо
неравенство
\begin{equation}\label{21a}J_n^{(\gamma)}\leqslant\gamma^{-\frac{n}{4}}\cdot\left(
\prod^n_{k=1} \alpha_k\right)^{1/2} \cdot \left[F
\left(\frac{2}{n}\sqrt{\gamma}\right)\right]^{\frac{n}{2}}.\end{equation}

Знак равенства в этом неравенстве достигается, когда $a_k$ и $B_k$,
$k=\overline{0,n}$, являются, соответственно, полюсами и круговыми
областями квадратичного дифференциала}
\begin{equation}\label{3a}Q(w)dw^2=-\frac{(n^2-\gamma)w^n+\gamma}{w^2(w^n-1)^2}\,dw^2.\end{equation}

\textbf{\textit{Доказательство теоремы 1.}} При доказательстве
теоремы 1 существенно используются идеи доказательства теоремы 1 из
работы \cite{8} и свойства разделяющего преобразования (см.
\cite{3}--\cite{6}, \cite{9}, \cite{10}). Положим
$$P_k:=P_k(A_n):=\{w\in\mathbb{C}\backslash\{0\}:\,\theta_k<\arg w<\theta_{k+1}\}.$$
Рассмотрим введеную ранее систему функций
$\zeta=\pi_k(w)=-i\left(e^{-i\theta_k}w\right)^\frac{1}{\alpha_k}$,\quad
$k=\overline{1,n}$. Пусть $\Omega_k^{(1)}$,
$k=\overline{1,n}$, обозначает область плоскости $\zeta$, полученную
в результате объединения связной компоненты множества
$\pi_k(B_k\bigcap\overline{P}_k)$, содержащей точку $\pi_k(a_k)$, со
своим симметричным отражением относительно мнимой оси. В свою
очередь, через $\Omega_k^{(2)}$, $k=\overline{1,n}$, обозначаем
область плоскости $\mathbb{C}_\zeta$, полученную в результате
объединения связной компоненты множества
$\pi_k(B_{k+1}\bigcap\overline{P}_k)$, содержащей точку
$\pi_k(a_{k+1})$, со своим симметричным отражением относительно мнимой
оси, $B_{n+1}:=B_1$, $\pi_n(a_{n+1}):=\pi_n(a_1)$. Кроме того,
$\Omega_k^{(0)}$ будет обозначать область плоскости
$\mathbb{C}_\zeta$ , полученную в результате объединения связной
компоненты множества $\pi_k(B_0\bigcap\overline{P}_k)$, содержащей
точку $\zeta=0$, со своим симметричным отражением относительно
мнимой оси. Обозначим $\pi_k(a_k):=\omega_k^{(1)}$,
$\pi_k(a_{k+1}):=\omega_k^{(2)}$, $k=\overline{1,n}$,
$\pi_n(a_{n+1}):=\omega_n^{(2)}$. Из определения функций $\pi_k$
вытекает, что
$$|\pi_k(w)-\omega_k^{(1)}|\sim\frac{1}{\alpha_k}|a_k|^{\frac{1}{\alpha_k}-1}\cdot|w-a_k|,\quad
w\rightarrow a_k,\quad w\in\overline{P_k},$$
$$|\pi_k(w)-\omega_k^{(2)}|\sim\frac{1}{\alpha_k}|a_{k+1}|^{\frac{1}{\alpha_k}-1}\cdot|w-a_{k+1}|,\quad
w\rightarrow a_{k+1},\quad w\in\overline{P_k},$$
$$|\pi_k(w)|\sim|w|^\frac{1}{\alpha_k},\quad
w\rightarrow 0,\quad w\in\overline{P_k}.$$ Тогда, используя
соответствующие результаты работ \cite{3}--\cite{6}, \cite{9},
получаем неравенства
\begin{equation}\label{4a}r\left(B_k,a_k\right)\leqslant\left[\frac{r\left(\Omega_k^{(1)},\omega_k^{(1)}\right)
\cdot
r\left(\Omega_k^{(2)},\omega_k^{(2)}\right)}{\frac{1}{\alpha_k}
|a_k|^{\frac{1}{\alpha_k}-1}\cdot\frac{1}{\alpha_{k-1}}
|a_k|^{\frac{1}{\alpha_{k-1}}-1}}\right]^\frac{1}{2},\end{equation}
$$k=\overline{1,n},\quad \Omega_0^{(2)}:=\Omega_n^{(2)},\quad \omega_0^{(2)}:=\omega_n^{(2)},$$
\begin{equation}\label{5a}r\left(B_0,0\right)\leqslant\left[\prod \limits_{k=1}^n
r^{\alpha_k^2}\left(\Omega_k^{(0)},0\right)\right]^\frac{1}{2}.\end{equation}
Повторяя рассуждения приведенные в \cite{9} при доказательстве
теоремы 5.2.1. с учетом введенных наборов областей
$\{P_{k}\}_{k=1}^{n}$, функций $\{\pi_{k}\}_{k=1}^{n}$ и чисел
$\{\theta_{k}\}_{k=1}^{n}$ получаем неравенство для исследуемого
функционала (\ref{1a})

$$J_n^{(\gamma)}\leqslant
\prod\limits_{k=1}^n\left[r\left(\Omega_k^{(0)},0\right)\right]^{\frac{\alpha_k^2}{2}\gamma}\cdot
\prod\limits_{k=1}^n\left[\frac{r\left(\Omega_{k-1}^{(2)},\omega_{k-1}^{(2)}\right)
r\left(\Omega_k^{(1)},\omega_k^{(1)}\right)}
{\frac{1}{\alpha_{k-1}\cdot\alpha_k}|a_k|^{\frac{1}{\alpha_{k-1}}-1}\cdot
|a_k|^{\frac{1}{\alpha_k}-1}}\right]^\frac{1}{2}=$$
\begin{equation}\label{6a}=\prod\limits_{k=1}^n\alpha_k\cdot \prod\limits_{k=1}^n
\frac{|a_k|}{|a_k a_{k+1}|^\frac{1}{2\alpha_k}}
\cdot\left[\prod\limits_{k=1}^n
r^{\gamma\alpha_k^2}\left(\Omega_k^{(0)},0\right)
\prod\limits_{k=1}^n r\left(\Omega_k^{(1)},\omega_k^{(1)}\right)
r\left(\Omega_k^{(2)},\omega_k^{(2)}\right)\right]^\frac{1}{2}.\end{equation}

Выражение, стоящее в скобках формулы (\ref{6a}),
представляет собой произведение значений функционала
$r^{\beta^2}(\Omega_k^{(0)},0)r(\Omega_k^{(1)},\omega_k^{(1)})r(\Omega_k^{(2)},\omega_k^{(2)})$
на тройках неналегающих областей
$\left(\Omega_k^{(0)},\Omega_k^{(1)},\Omega_k^{(2)}\right)$
плоскости $\zeta$.

Известно \cite{7}, что функционал
$$Y_3(\sigma_1,\sigma_2,\sigma_3)=
\frac{r^{\sigma_1}(D_1,d_1)\cdot r^{\sigma_2}(D_2,d_2)\cdot
r^{\sigma_3}(D_3,d_3)}
{|d_1-d_2|^{\sigma_1+\sigma_2-\sigma_3}\cdot|d_1-d_3|^{\sigma_1-\sigma_2+\sigma_3}\cdot
|d_2-d_3|^{-\sigma_1+\sigma_2+\sigma_3}},$$
$\sigma_k\in\mathbb{R}^+$,\, $d_k\in
D_k\subset\overline{\mathbb{C}}$,\, $D_k\bigcap D_p=\varnothing$,\,
$k=1,2,3$,\, $p=1,2,3$,\, $k\neq p$, инвариантен относительно всех
конформных автоморфизмов комплексной плоскости
$\overline{\mathbb{C}}$.

С учетом этого справедлива следующая оценка
$$J_n^{(\gamma)}\leqslant
\left(\prod\limits_{k=1}^n\alpha_k\right)\cdot\prod\limits_{k=1}^n
\frac{|a_k|}{|a_k a_{k+1}|^\frac{1}{2\alpha_k}}\times$$
\begin{equation}\label{73a}\times \left\{\prod\limits_{k=1}^n\frac{
r^{\gamma\alpha_k^2}\left(\Omega_k^{(0)},0\right)\cdot
r\left(\Omega_k^{(1)},\omega_k^{(1)}\right)\cdot
r\left(\Omega_k^{(2)},\omega_k^{(2)}\right)}{|\omega_k^{(1)}\cdot
\omega_k^{(2)}|^{\gamma\alpha_k^2}|\omega_k^{(1)}-
\omega_k^{(2)}|^{2-\gamma\alpha_k^2}}\right\}^\frac{1}{2}\times\end{equation}
$$\times\left[\prod\limits_{k=1}^n |\omega_k^{(1)}\cdot
\omega_k^{(2)}|^{\gamma\alpha_k^2}|\omega_k^{(1)}-
\omega_k^{(2)}|^{2-\gamma\alpha_k^2}\right]^\frac{1}{2}.$$

Справедливы соотношения $|\omega_k^{(1)}|=|a_k|^\frac{1}{\alpha_k}$,
$|\omega_k^{(2)}|=|a_{k+1}|^\frac{1}{\alpha_k}$,
$|\omega_k^{(1)}-\omega_k^{(2)}|=|a_k|^\frac{1}{\alpha_k}+|a_{k+1}|^\frac{1}{\alpha_k}$.

Далее, учитывая эти равенства имеем

$$J_n^{(\gamma)}\leqslant \left(\prod\limits_{k=1}^n\alpha_k\right)\cdot\prod\limits_{k=1}^n
\frac{|a_k|}{|a_k a_{k+1}|^\frac{1}{2\alpha_k}}\times$$
$$\times\left(\prod\limits_{k=1}^n|\omega_k^{(1)}-
\omega_k^{(2)}|\right)\left(\prod\limits_{k=1}^n\frac{|\omega_k^{(1)}\cdot
\omega_k^{(2)|}}{|\omega_k^{(1)}-
\omega_k^{(2)}|}\right)^{\frac{\gamma\alpha_{k}^{2}}{2}}\times$$
$$\times \left\{\prod\limits_{k=1}^n\frac{
r^{\gamma\alpha_k^2}\left(\Omega_k^{(0)},0\right)\cdot
r\left(\Omega_k^{(1)},\omega_k^{(1)}\right)\cdot
r\left(\Omega_k^{(2)},\omega_k^{(2)}\right)}{|\omega_k^{(1)}\cdot
\omega_k^{(2)}|^{\gamma\alpha_k^2}|\omega_k^{(1)}-
\omega_k^{(2)}|^{2-\gamma\alpha_k^2}}\right\}^\frac{1}{2}=$$
\begin{equation}\label{7a}\end{equation}
$$=2^{n}\cdot\left(\prod\limits_{k=1}^n\alpha_k\right)\cdot\prod\limits_{k=1}^n
\chi\left(\left|\frac{a_k}{a_{k+1}}\right|^\frac{1}{2\alpha_k}\right)|a_k|\times$$
$$\times2^{-\frac{\gamma}{2}\sum\limits_{k=1}^n\alpha_k}\left[\prod\limits_{k=1}^n
\chi\left(\left|\frac{a_k}{a_{k+1}}\right|^\frac{1}{2\alpha_k}\right)\right]^{-\frac{\gamma\alpha_{k}^{2}}{2}}
\left(\prod\limits_{k=1}^n\left|\frac{a_{k+1}}{a_{k}}\right|\right)^{\frac{\gamma\alpha_{k}^{2}}{2}}\times$$
$$\times \left\{\prod\limits_{k=1}^n\frac{
r^{\gamma\alpha_k^2}\left(\Omega_k^{(0)},0\right)\cdot
r\left(\Omega_k^{(1)},\omega_k^{(1)}\right)\cdot
r\left(\Omega_k^{(2)},\omega_k^{(2)}\right)}{|\omega_k^{(1)}\cdot
\omega_k^{(2)}|^{\gamma\alpha_k^2}|\omega_k^{(1)}-
\omega_k^{(2)}|^{2-\gamma\alpha_k^2}}\right\}^\frac{1}{2}=$$
$$=2^{n-\frac{\gamma}{2}\sum\limits_{k=1}^n\alpha_k^2}\cdot
\left(\prod\limits_{k=1}^n\alpha_k\right)
\cdot\prod\limits_{k=1}^n\left[\chi\left(\Bigl|\frac{a_k}{a_{k+1}}\Bigr|
^\frac{1}{2\alpha_k}\right)\right]^{1-\frac{\gamma\alpha_k^2}{2}}\times$$
$$\times\prod\limits_{k=1}^n|a_k|^{1+\frac{1}{4}\gamma(\alpha_k+\alpha_{k-1})}\times$$
$$\times\left\{\prod\limits_{k=1}^n\frac{
r^{\gamma\alpha_k^2}\left(\Omega_k^{(0)},0\right)\cdot
r\left(\Omega_k^{(1)},\omega_k^{(1)}\right)\cdot
r\left(\Omega_k^{(2)},\omega_k^{(2)}\right)}{|\omega_k^{(1)}\cdot
\omega_k^{(2)}|^{\gamma\alpha_k^2}|\omega_k^{(1)}-
\omega_k^{(2)}|^{2-\gamma\alpha_k^2}}\right\}^\frac{1}{2}=$$

$$=2^{n-\frac{\gamma}{2}\sum\limits_{k=1}^n\alpha_k^2}\cdot
\left(\prod\limits_{k=1}^n\alpha_k\right)\cdot\mathcal{L}^{(\gamma)}\left(A_n\right)\times$$
$$\times\left\{\prod\limits_{k=1}^n\frac{
r^{\gamma\alpha_k^2}\left(\Omega_k^{(0)},0\right)\cdot
r\left(\Omega_k^{(1)},\omega_k^{(1)}\right)\cdot
r\left(\Omega_k^{(2)},\omega_k^{(2)}\right)}{|\omega_k^{(1)}\cdot
\omega_k^{(2)}|^{\gamma\alpha_k^2}|\omega_k^{(1)}-
\omega_k^{(2)}|^{2-\gamma\alpha_k^2}}\right\}^\frac{1}{2}.$$

При каждом $k=\overline{1,n}$ несложно указать конформный
автоморфизм $\zeta=T_k(z)$ плоскости комплексных чисел
$\overline{\mathbb{C}}$ такой, что $T_k(0)=0$,\,
$T_k\left(\omega_k^{(s)}\right)=(-1)^s\cdot i$,\,
$G_k^{(q)}:=T_k\left(\Omega_k^{(q)}\right)$,\, $k=\overline{1,n}$,\,
$s=1,2$,\, $q=0,1,2$. Тогда

$$\left\{\prod\limits_{k=1}^n\frac{
r^{\gamma\alpha_k^2}\left(\Omega_k^{(0)},0\right)\cdot
r\left(\Omega_k^{(1)},\omega_k^{(1)}\right)\cdot
r\left(\Omega_k^{(2)},\omega_k^{(2)}\right)}{|\omega_k^{(1)}\cdot
\omega_k^{(2)}|^{\gamma\alpha_k^2}|\omega_k^{(1)}-
\omega_k^{(2)}|^{2-\gamma\alpha_k^2}}\right\}^\frac{1}{2}=$$

$$=\left\{\prod\limits_{k=1}^n\frac{r^{\alpha_k^2\gamma}\left(G_k^{(0)},0\right)\cdot
r\left(G_k^{(1)},-i\right)\cdot
r\left(G_k^{(2)},i\right)}{2^{2-\gamma\alpha_{k}^{2}}}\right\}^\frac{1}{2}.$$

Далее, иcпользуя результаты работ \cite{7} -- \cite{9} получим следующее выражение
$$J_n^{(\gamma)}\leqslant
2^{n-\frac{\gamma}{2}\sum\limits_{k=1}^n\alpha_k^2}\cdot
\left(\prod\limits_{k=1}^n\alpha_k\right)\cdot\mathcal{L}^{(\gamma)}\left(A_n\right)\times$$

$$\times\prod\limits_{k=1}^n\left\{\frac{r^{\alpha_k^2\gamma}\left(G_k^{(0)},0\right)\cdot
r\left(G_k^{(1)},-i\right)\cdot
r\left(G_k^{(2)},i\right)}{2^{2-\gamma\alpha_{k}^{2}}}\right\}^\frac{1}{2}=$$
$$=2^{n-\frac{\gamma}{2}\sum\limits_{k=1}^n\alpha_k^2}\left(\prod\limits_{k=1}^n\alpha_k\right)\cdot\mathcal{L}^{(\gamma)}(A_n)\cdot2^{-n+\frac{\gamma}{2}\sum\limits_{k=1}^n\alpha_k^2}\times$$
$$\times\left[\prod\limits_{k=1}^n
r^{\alpha_k^2\gamma}\left(G_k^{(0)},0\right)\cdot
r\left(G_k^{(1)},-i\right)\cdot
r\left(G_k^{(2)},i\right)\right]^\frac{1}{2}\leqslant$$
\begin{equation}\label{100a}\leqslant\left(\prod\limits_{k=1}^n\alpha_k\right)\cdot\mathcal{L}^{(\gamma)}(A_n)\cdot
\left[\prod\limits_{k=1}^n
r^{\alpha_k^2\gamma}\left(G_k^{(0)},0\right)\cdot
r\left(G_k^{(1)},-i\right)\cdot
r\left(G_k^{(2)},i\right)\right]^\frac{1}{2}.\end{equation}

В результате проведенных вычислений получим неравенство для
функционала (1)
\begin{equation}\label{122a}J_n^{(\gamma)}\leqslant\gamma^{-\frac{n}{4}}\left(
\prod\limits^n_{k=1} \alpha_k\right)^{1/2}
\left[\prod\limits^n_{k=1}
F\left(\frac{2}{n}\sqrt{\gamma}\right)\right]^{1/2}.\end{equation}

Теперь рассмотрим вспомогательный функционал $$\tilde
J_n^{(\gamma)}= \left(\prod\limits^n_{k=1}
\alpha_k\right)^{-1/2}\cdot J_n^{(\gamma)}.$$ Из неравенства
(\ref{122a}) следует, что
$$\tilde J_n= \gamma^{-\frac{n}{4}}\left[\prod\limits^n_{k=1}F\left(\frac{2}{n}\sqrt{\gamma}\right)\right] ^{\frac{1}{2}}
.$$

Аналогично \cite{4}, рассмотрим экстремальную задачу

\begin{equation}\label{10a}\prod^n_{k=1}F(\alpha_k) \longrightarrow \max; \quad
\sum^n_{k=1}\alpha_k=2.\end{equation}

Ясно, что необходимые условия экстремума имеют вид
\begin{equation}\label{11a}\frac{F'(\alpha_k)}{F(\alpha_k)}=-\frac{\lambda}
{\prod\limits^n_{k=1}F(\alpha_k)}, \quad k=\overline{1,n}.
\end{equation} ($\lambda$ -- фиксированное, вещественное число). Функция
$\Phi(\alpha)= \frac{F'(\alpha)}{F(\alpha)}$ убывает на промежутке
$(0, t_0], 1,168 < t_0 < 1,170$ и возрастает на $[t_0, 2)$. По аналогии
с работой \cite{4} рассмотрим функцию
$H(\alpha)=\Phi(\alpha)-\Phi(2-\alpha)$, $\alpha \in (0,2)$.
Несложные вычисления показывают, что функция $H(\alpha)$
положительна на интервале $(0,1)$. Обозначим
$\alpha_0=\max\limits_{1\leq k\leq n}\alpha_k$. Пусть
$\alpha_0=\alpha_{k_0}, \ 1 \leq k_0 \leq n, \ k_0 \in \mathbb{N}$.
Допустим, что $\alpha_0 > 1$, тогда
$\alpha_k \leq 2-\alpha_0 < 1$ для всех $k=\overline{1,n}, \\
k\not = k_0$. Следуя работе \cite{4} получим соотношения $$
\Phi(\alpha_0)=\Phi(2-(2-\alpha_0)) < \Phi(2-\alpha_0)
 \leq \Phi(\alpha_k), k=\overline{1,n},  k\not=0.$$

Последнее неравенство противоречит необходимим условиям (\ref{11a}).
Поэтому все $\alpha_k \in (0,1), \ k=\overline{1,n}$.

Тогда приходим к выводу $$\tilde J_n^{(\gamma)}\leq
\gamma^{-\frac{n}{4}}\left[ F
\left(\frac{2}{n}\sqrt{\gamma}\right)\right]^{n/2}.$$

Возвращаясь к исходному функционалу (1) получаем окончательное
неравенство $$J_n^{(\gamma)} \leq\gamma^{-\frac{n}{4}}
\left[\prod^n_{k=1}\alpha_k\right]^{1/2}\left[F
\left(\frac{2}{n}\sqrt{\gamma}\right)\right]^{\frac{n}{2}}.$$

Утверждение о знаке равенства
проверяется непосредственно.

Из теоремы 1 непосредственно вытекает теорема 4 работы \cite{4}.

\textbf{Следствие 1.} Пусть $n \in \mathbb{N}$, $n\geq 2$. Тогда для
любой системы различных точек единичной окружности
$A_n=\{a_k\}^n_{k=1}$ и любого набора взаимно непересекающихся
областей  $\{B_k\}^n_{k=0}$, $0 \in B_0$, $a_k \in B_k$,
$k=\overline{1,n}$ справедливо неравенство

$$r(B_0,0)\prod^n_{k=1}r(B_k, a_k)\leqslant\frac{4^{n+\frac{\gamma}{n}}
\gamma^\frac{\gamma}{n}n^n}{(n^2-\gamma)^{n+\frac{\gamma}{n}}}
\left(\frac{n-\sqrt{\gamma}}{n+\sqrt{\gamma}}\right)^{2\sqrt{\gamma}}.$$
Знак равенства достигается при условиях теоремы 1.

При $\gamma=1$ и $n\geq 2$ следствие 1 было получено В. Н. Дубининым
в работе (\cite{4}, 1988 г.), причем из его метода следует, что результат верен и при $0<\gamma\leq1$ Позднее Г.В.
Кузьмина повторила этот результат для односвязных областей другим
методом.

Из неравенства (\ref{100a}) и метода доказательства теоремы 1
получаем следующие утверждения.

\textbf{Следствие 2.} \cite{13} Пусть $n\in\mathbb{N}$, $n\geq 2$, $\gamma\in
(0,1]$. Тогда для любой системы различных точек единичной окружности
$A_n=\{a_k\}^n_{k=1}$ и любого набора взаимно непересекающихся
областей  $\{B_k\}^n_{k=0}$, $0 \in B_0$, $a_k \in B_k$,
$k=\overline{1,n}$ справедливо неравенство

$$J_n^{(\gamma)}\leq\left( \prod^n_{k=1}
\alpha_k\right)^{1/2} \cdot 2^n
\left(\frac{2}{n}\right)^{\frac{n}{2}}
\frac{\left(\frac{4\gamma}{n^2}\right)^{\frac{\gamma}{n}}}
{\left(1-\frac{\gamma}{n^2}\right)^{n+\frac{\gamma}{n}}}
\left(\frac{1-\frac{\sqrt{\gamma}}{n}}{1+\frac{\sqrt{\gamma}}{n}}
\right)^{2\sqrt{\gamma}}.$$ Знак равенства достигается при условиях
теоремы 1.

\textbf{Следствие 3.} \cite{8} Пусть $n\in \mathbb{N}$, $n\geqslant 2$ и
$\gamma\in (0; 0,2]$. Тогда при условиях теоремы 1 справедливо
неравенство $$J_n^{(\gamma)}\leqslant\left( \prod^n_{k=1}
\alpha_k\right)\cdot 2^n \left(\frac{2}{n}\right)^{\frac{n}{2}}
\frac{\left(\frac{4\gamma}{n^2}\right)^{\frac{\gamma}{n}}}
{\left(1-\frac{\gamma}{n^2}\right)^{n+\frac{\gamma}{n}}}
\left(\frac{1-\frac{\sqrt{\gamma}}{n}}{1+\frac{\sqrt{\gamma}}{n}}
\right)^{2\sqrt{\gamma}}.$$ Знак равенства достигается при условиях
теоремы 1.

\textbf{Следствие 4.} Пусть $n\in \mathbb{N}$, $n\geqslant 2$ и
$\gamma\in (0, 1]$. Тогда для любой $n$-лучевой системы точек
$A_n=\{a_k\}_{k=1}^n$ и любого набора взаимно непересекающихся
областей $B_k$, $a_k\in B_k\subset\overline{\mathbb{C}}$,
$k=\overline{0,n}$, справедливо неравенство
$$J_n^{(\gamma)}\leqslant\frac{4^{n+\frac{\gamma}{n}}
\gamma^\frac{\gamma}{n}n^n}{(n^2-\gamma)^{n+\frac{\gamma}{n}}}
\left(\frac{n-\sqrt{\gamma}}{n+\sqrt{\gamma}}\right)^{2\sqrt{\gamma}}.$$
Знак равенства в этом неравенстве достигается, когда $a_k$ и $B_k$,
$k=\overline{0,n}$, являются, соответственно, полюсами и круговыми
областями квадратичного дифференциала
$$Q(w)dw^2=-\frac{(n^2-\gamma)w^n+R^{n}\gamma}{w^2(w^n-R^{n})^2}\,dw^2,$$
где $R^{n+\gamma}=\mathcal{L}^{(\gamma)}\left(A_n\right)$.

\begin{center} {\bf ЛИТЕРАТУРА}
\end{center}

\end{document}